\documentclass[a4paper,11pt]{amsart}
\usepackage{graphicx}
 \usepackage{mathptmx}
\usepackage{amsmath}
\usepackage{amssymb}
\usepackage{enumitem}
\usepackage{xcolor}

\newtheorem{theorem}{Theorem}[section]

\newtheorem{corollary}[theorem]{Corollary}
\newtheorem{proposition}[theorem]{Proposition}

\begin{document}

\title[A note on degenerate Euler and Bernoulli polynomials of complex variable]{A note on degenerate Euler and Bernoulli polynomials of complex variable}

\author{Dae San  Kim }
\address{Department of Mathematics, Sogang University, Seoul 121-742, Republic of Korea}
\email{dskim@sogang.ac.kr}
\author{Taekyun  Kim}
\address{Department of Mathematics, Kwangwoon University, Seoul 139-701, Republic of Korea}
\email{tkkim@kw.ac.kr}
\author{Hyunseok Lee}
\address{Department of Mathematics, Kwangwoon University, Seoul 139-701, Republic of Korea}
\email{luciasconstant@gmail.com}


\maketitle

\begin{abstract}
Recently, the so called the new type Euler polynomials have been studied without considering Euler polynomials of complex variable. Here we study the degenerate versions of these new type Euler polynomials, namely degenerate cosine-Euler polynomials and degenerate sine-Euler polynomials and also the corresponding ones for Bernoulli polynomials, namely degenerate cosine-Bernoulli polynomials and degenerate sine-Bernoulli polynomials by considering the degenerate Euler polynomials of complex variable and the degenerate Bernoulli polynomials of complex variable. We derived some properties and identities for those new polynomials. Here we note that our result gives an affirmative answer to the question raised by the reviewer of the paper [15].
\end{abstract}

\section{Introduction} 
As is well known, the ordinary Bernoulli polynomials $B_{n}(x)$ and Euler polynomials $E_{n}(x)$ are respectively defined by 
\begin{equation}\label{1}
	\frac{t}{e^{t}-1}e^{xt}=\sum_{n=0}^{\infty}B_{n}(x)\frac{t^{n}}{n!}, 
\end{equation}
and 
\begin{equation}\label{2}
	\frac{2}{e^{t}+1}e^{xt}=\sum_{n=0}^{\infty}E_{n}(x) \frac{t^{n}}{n!} ,\qquad \textrm{(see [1-20]).}
\end{equation}
For any nonzero $\lambda\in\mathbb{R}$, the degenerate exponential function is defined by 
\begin{equation}\label{3}
	e_{\lambda}^{x}(t)=(1+\lambda t)^{\frac{x}{\lambda}},\quad e_{\lambda}(t)=e_{\lambda}^{1}(t),\quad (\mathrm{see}\ [8]).
\end{equation}
In [1,2], Carlitz considered the degenerate Bernoulli and Euler polynomials which are given by 
\begin{equation}\label{4}
	\frac{t}{e_{\lambda}(t)-1}e_{\lambda}^{x}(t)=\frac{t}{(1+\lambda t)^{\frac{1}{\lambda}}-1}(1+\lambda t)^{\frac{x}{\lambda}}=\sum_{n=0}^{\infty}\beta_{n,\lambda}(x) \frac{t^{n}}{n!}
\end{equation}
and 
\begin{equation}\label{5}
	\frac{2}{e_{\lambda}(t)+1}e_{\lambda}^{x}(t)=\frac{2}{(1+\lambda t)^{\frac{1}{\lambda}}+1}(1+\lambda t)^{\frac{x}{\lambda}}=\sum_{n=0}^{\infty}\mathcal{E}_{n,\lambda}(x) \frac{t^{n}}{n!}.
\end{equation}
Note that 
\begin{displaymath}
	\lim_{\lambda\rightarrow 0}\beta_{n,\lambda}(x)=B_{n}(x),\quad \lim_{\lambda\rightarrow 0}\mathcal{E}_{n,\lambda}(x)=E_{n}(x).
\end{displaymath}
The falling factorial sequence is defined as 
\begin{displaymath}
	(x)_{0}=1,\quad (x)_{n}=x(x-1)\cdots(x-n+1),\ (n\ge 1),\quad (\mathrm{see\ [17]}).
\end{displaymath}
The Stirling numbers of the first kind are defined by the coefficients in the expansion of $(x)_n$ in terms of powers of $x$ as follows:
\begin{equation}\label{6}
	(x)_{n}=\sum_{l=0}^{n}S^{(1)}(n,l)x^{l},\quad (\mathrm{see},\ [7,11,17]).
\end{equation}
The Stirling numbers of the second kind are defined by 
\begin{equation}\label{7}
	x^{n}=\sum_{l=0}^{n}S^{(2)}(n,l)(x)_{l},\ (n\ge 0),\ (\mathrm{see}\ [9,10,17])
.\end{equation}
In [9], the degenerate stirling numbers of the second kind are defined by the generating function 
\begin{equation}\label{8}
	\frac{1}{k!}\big(e_{\lambda}(t)-1\big)^{k}=\sum_{n=k}^{\infty}S_{\lambda}^{(2)}(n,k)\frac{t^{n}}{n!},\ (k\ge 0). 
\end{equation}
Note that $\displaystyle\lim_{\lambda\rightarrow 0}S_{\lambda}^{(2)}(n,k)=S^{(2)}(n,k),\ (n,k\ge 0)\displaystyle$. \\
~~\\
Recently, Masjed-Jamei, Beyki and Koepf introduced the new type Euler polynomials which are given by
\begin{equation}\label{9}
	\frac{2e^{pt}}{e^{t}+1}\cos qt=\sum_{n=0}^{\infty}E_{n}^{(c)}(p,q)\frac{t^{n}}{n!}, 
\end{equation}
\begin{equation}\label{10}
	\frac{2e^{pt}}{e^{t}+1}\sin qt=\sum_{n=0}^{\infty}E_{n}^{(s)}(p,q)\frac{t^{n}}{n!},\quad (\mathrm{see}\ [15]).
\end{equation}
They also considered the cosine-polynomials and sine-polynomials defined by
\begin{equation}\label{11}
	e^{pt}\cos qt=\sum_{n=0}^{\infty}C_{n}(p,q)\frac{t^{n}}{n!}, 
\end{equation}
and 
\begin{equation}\label{12}
	e^{pt}\sin qt=\sum_{n=0}^{\infty}S_{n}(p,q)\frac{t^{n}}{n!},\quad (\mathrm{see}\ [15]).
\end{equation}
In [15], the authors deduced many interesting identities and properties for those polynomials. \\
~~\\
It is well known that 
\begin{equation}\label{13}
	e^{ix}=\cos x+i\sin x,\quad\mathrm{where}\,\,x \in \mathbb{R},\ i=\sqrt{-1},\quad (\mathrm{see}\ [20]).
\end{equation}
From \eqref{1} and \eqref{2}, we note that 
\begin{equation}\label{14}
\frac{t}{e^{t}-1}e^{(x+iy)t}=\sum_{n=0}^{\infty}B_{n}(x+iy)\frac{t^{n}}{n!}, 
\end{equation}
and 
\begin{equation}\label{15}
\frac{2}{e^{t}+1}e^{(x+iy)t}=\sum_{n=0}^{\infty}E_{n}(x+iy)\frac{t^{n}}{n!}.
\end{equation}
By \eqref{14} and \eqref{15}, we get 
\begin{align}\label{16}
\frac{t}{e^{t}-1}e^{xt}\cos yt &=\sum_{n=0}^{\infty}\frac{B_{n}(x+iy)+B_{n}(x-iy)}{2}\frac{t^{n}}{n!}=\sum_{n=0}^{\infty}B_n^{(c)}(x,y)\frac{t^n}{n!}, \\
\frac{t}{e^{t}-1}e^{xt}\sin yt &=\sum_{n=0}^{\infty}\frac{B_{n}(x+iy)-B_{n}(x-iy)}{2i}\frac{t^{n}}{n!}=\sum_{n=0}^{\infty}B_n^{(s)}(x,y)\frac{t^n}{n!},\nonumber\\
\frac{2}{e^{t}+1}e^{xt}\cos yt &=\sum_{n=0}^{\infty}\frac{E_{n}(x+iy)+E_{n}(x-iy)}{2}\frac{t^{n}}{n!}=\sum_{n=0}^{\infty}E_n^{(c)}(x,y)\frac{t^n}{n!},\nonumber
\end{align}
and 
\begin{displaymath}
\quad\frac{2}{e^{t}+1}e^{xt}\sin yt =\sum_{n=0}^{\infty}\frac{E_{n}(x+iy)-E_{n}(x-iy)}{2i}\frac{t^{n}}{n!}=\sum_{n=0}^{\infty}E_n^{(s)}(x,y)\frac{t^n}{n!},\,[\mathrm{see}\ 12].
\end{displaymath}

\vspace{0.1in}

In view of \eqref{4} and \eqref{5}, we study the degenerate Bernoulli and Euler polynomials with complex variable and investigate some identities and properties for those polynomials. The outline of this paper is as follows. In Section 1, we will beriefly recall the degenerate Bernoulli and Euler polynomials of Carlitz and the degenerate Stirling numbers of the second kind. Then we will introduce so called the new type Euler polynomials, and the cosine-polynomials and sine-polynomials recently introduced in [15]. Then we indicate that the new type Euler polynomials and the corresponding Bernoulli polynomials can be expressed by considering Euler and Bernoulli polynomials of complex variable and treating the real and imaginary parts separately. In Section 2, the degenerate cosine-polynomials and degenerate sine-polynomials were introduced and their explicit expressions were derived. The degenerate cosine-Euler polynomials and degenerate sine-Euler polynomials were expressed in terms of degenerate cosine-polynomials and degenerate sine-polynomials and vice versa. Further, some reflection identities were found for the degenerate cosine-Euler polynomials and degenerate sine-Euler polynomials. In Section 3, the degenerate cosine-Bernoulli polynomials and degenerate sine-Bernoulli polynomials were introduced. They were expressed in terms of degenerate cosine-polynomials and degenerate sine-polynomials and vice versa. Reflection symmetries were deduced for the degenerate cosine-Bernoulli polynomials and degenrate sine-Bernoulli polynomials.

\section{Degenerate Euler polynomials of complex variable} 

Here we will consider the degenerate Euler polynomials of complex variable and, by treating the real and imaginary parts separately, introduce the degenerate cosine-Euler polynomials and degenerate sine-Euler polynomials. They are degenerate versions of the new type Euler polynomials studied in [15].

The degenerate sine and cosine functions are defined by 
\begin{equation}\label{17}
	\cos_{\lambda}t=\frac{e_{\lambda}^{i}(t)+e_{\lambda}^{-i}(t)}{2},\quad \sin_{\lambda}t= \frac{e_{\lambda}^{i}(t)-e_{\lambda}^{-i}(t)}{2i}.
\end{equation}
From \eqref{13}, we note that 
\begin{displaymath}
	\lim_{\lambda\rightarrow 0}\cos_{\lambda}t=\cos t,\quad \lim_{\lambda\rightarrow 0}\sin_{\lambda}t=\sin t.  
\end{displaymath}
By \eqref{5}, we get 
\begin{equation}\label{18}
	\frac{2}{e_{\lambda}(t)+1}e_{\lambda}^{x+iy}(t)=\sum_{n=0}^{\infty}\mathcal{E}_{n,\lambda}(x+iy)\frac{t^{n}}{n!}.
\end{equation}
and
\begin{equation}\label{19}
	\frac{2}{e_{\lambda}(t)+1}e_{\lambda}^{x-iy}(t)=\sum_{n=0}^{\infty}\mathcal{E}_{n,\lambda}(x-iy)\frac{t^{n}}{n!}.
\end{equation}
Now, we define the degenerate cosine and degenerate sine function as
\begin{equation}\label{20}
\cos_{\lambda}^{(y)}(t)=\frac{e_{\lambda}^{iy}(t) + e_{\lambda}^{-iy} (t)}{2}=\cos\bigg(\frac{y}{\lambda}\log(1+\lambda t)\bigg),
\end{equation}

\begin{equation}\label{21}
\sin_{\lambda}^{(y)}(t)=\frac{e_{\lambda}^{iy}(t) - e_{\lambda}^{-iy} (t)}{2i}=\sin\bigg(\frac{y}{\lambda}\log(1+\lambda t)\bigg). 
\end{equation}
Note that $\displaystyle\lim_{\lambda\rightarrow 0}\cos_{\lambda}^{(y)}(t)=\cos yt,\ \lim_{\lambda\rightarrow 0}\sin_{\lambda}^{(y)}(t)=\sin yt\displaystyle$. \\
\noindent From \eqref{18} and \eqref{19}, we note that 
\begin{equation}\label{22}
	\frac{2}{e_{\lambda}(t)+1}e_{\lambda}^{x}(t)\cos_{\lambda}^{(y)}(t)=\sum_{n=0}^{\infty}\bigg(\frac{\mathcal{E}_{n,\lambda}(x+iy)+\mathcal{E}_{n,\lambda}(x-iy)}{2}\bigg)\frac{t^{n}}{n!}, 
\end{equation}
and 
\begin{equation}\label{23}
	\frac{2}{e_{\lambda}(t)+1}e_{\lambda}^{x}(t)\sin_{\lambda}^{(y)}(t)=\sum_{n=0}^{\infty}\bigg(\frac{\mathcal{E}_{n,\lambda}(x+iy)-\mathcal{E}_{n,\lambda}(x-iy)}{2i}\bigg)\frac{t^{n}}{n!}.
\end{equation}
In view of \eqref{9} and \eqref{10}, we define the degenerate cosine-Euler polynomials and degenerate sine-Euler polynomials respectively by 
\begin{equation}\label{24}
	\frac{2}{e_{\lambda}(t)+1}e_{\lambda}^{x}(t)\cos_{\lambda}^{(y)}(t)=\sum_{n=0}^{\infty}\mathcal{E}_{n,\lambda}^{(c)}(x,y)\frac{t^{n}}{n!}, 
\end{equation}
and 
\begin{equation}\label{25}
	\frac{2}{e_{\lambda}(t)+1}e_{\lambda}^{x}(t)\sin_{\lambda}^{(y)}(t)=\sum_{n=0}^{\infty}\mathcal{E}_{n,\lambda}^{(s)}(x,y)\frac{t^{n}}{n!}.
\end{equation}
Note that $\displaystyle\lim_{\lambda\rightarrow 0}\mathcal{E}_{n,\lambda}^{(c)}(x,y)=E_{n}^{(c)}(x,y),\quad\lim_{\lambda\rightarrow 0}\mathcal{E}_{n,\lambda}^{(s)}(x,y)=E_{n}^{(s)}(x,y),\ (n\ge 0)\displaystyle$, where $E_{n}^{(c)}(x,y)$ and $E_{n}^{(s)}(x,y)$ are the new type of Euler polynomials of Masjed-Jamei, Beyki and Koepf (see [15]). \\
~~~\\
From \eqref{22}-\eqref{25}, we note that 
\begin{equation}\label{26}
	\mathcal{E}_{n,\lambda}^{(c)}(x,y)=\frac{\mathcal{E}_{n,\lambda}(x+iy)+\mathcal{E}_{n,\lambda}(x-iy)}{2}, 
\end{equation}
and 
\begin{equation}\label{27}
	\mathcal{E}_{n,\lambda}^{(s)}(x,y)=\frac{\mathcal{E}_{n,\lambda}(x+iy)-\mathcal{E}_{n,\lambda}(x-iy)}{2i},\quad (n\ge 0). 
\end{equation}
We recall here that the generalized falling factorial sequence is defined by
\begin{displaymath}
	(x)_{0,\lambda}=1,\quad (x)_{n,\lambda}=x(x-\lambda)(x-2\lambda)\cdots(x-(n-1)\lambda),\quad (n\ge 1). 
\end{displaymath}
Note that $\displaystyle\lim_{\lambda\rightarrow 1}(x)_{n,\lambda}=(x)_{n},\quad \lim_{\lambda\rightarrow 0}(x)_{n,\lambda}=x^{n}\displaystyle$.\\
We observe that  
\begin{align}\label{28}
e_{\lambda}^{iy}(t)&=(1+\lambda t)^{\frac{iy}{\lambda}}=e^{\frac{iy}{\lambda}\log(1+\lambda t)}\\
&=\sum_{k=0}^{\infty}\bigg(\frac{iy}{\lambda}\bigg)^{k}\frac{1}{k!}\big(\log(1+\lambda t)\big)^{k}\nonumber\\
&=\sum_{k=0}^{\infty}\lambda^{-k}(iy)^{k}\sum_{n=k}^{\infty}S^{(1)}(n,k)\frac{\lambda^{n}}{n!}t^{n}\nonumber\\
&=\sum_{n=0}^{\infty}\bigg(\sum_{k=0}^{n}\lambda^{n-k}i^{k}y^{k}S^{(1)}(n,k)\bigg)\frac{t^{n}}{n!}.\nonumber
\end{align}
From \eqref{20}, we can derive the following equation. 
\begin{align}\label{29}
	\cos_{\lambda}^{(y)}(t)&=\frac{e_{\lambda}^{iy}(t)+e_{\lambda}^{-iy}(t)}{2}\\
	&=\frac{1}{2}\sum_{n=0}^{\infty}\bigg(\sum_{k=0}^{n}\lambda^{n-k}(i^{k}+(-i)^{k})y^{k}S^{(1)}(n,k)\bigg)\frac{t^{n}}{n!}\nonumber\\
	&=\sum_{n=0}^{\infty}\bigg(
	\sum_{k=0}^{[\frac{n}{2}]}\lambda^{n-2k}(-1)^{k}y^{2k}S^{(1)}(n,2k)\bigg)\frac{t^{n}}{n!}\nonumber\\
	&=\sum_{k=0}^{\infty}\bigg(\sum_{n=2k}^{\infty}\lambda^{n-2k}(-1)^{k}y^{2k}S^{(1)}(n,2k)\bigg)\frac{t^{n}}{n!}.\nonumber
\end{align}
Note that 
\begin{displaymath}
	\lim_{\lambda\rightarrow 0}\cos_{\lambda}^{(y)}(t)=\sum_{k=0}^{\infty}(-1)^{k}y^{2k}\frac{t^{2k}}{(2k)!}=\cos yt. 
\end{displaymath}
By \eqref{21}, we get 
\begin{align}\label{30}	
\sin_{\lambda}^{(y)}(t)&=\frac{e_{\lambda}^{iy}(t)-e_{\lambda}^{-iy}(t)}{2i}\\
&=\frac{1}{2i}\sum_{n=0}^{\infty}\bigg(\sum_{k=0}^{n}\lambda^{n-k}(i^{k}-(-i)^{k})y^{k}S^{(1)}(n,k)\bigg)\frac{t^{n}}{n!}\nonumber\\ 
&=\sum_{n=1}^{\infty} \bigg(\sum_{k=0}^{[\frac{n-1}{2}]}\lambda^{n-2k-1}(-1)^{k}y^{2k+1}S^{(1)}(n,2k+1)\bigg)\frac{t^{n}}{n!} \nonumber	\\
&=\sum_{k=0}^{\infty}\bigg(\sum_{n=2k+1}^{\infty}(-1)^{k}\lambda^{n-2k-1}S^{(1)}(n,2k+1)\frac{t^{n}}{n!}
\bigg)y^{2k+1},\nonumber
\end{align}
where $[x]$ denotes the greatest integer $\leq x$. \\
~~\\
Note that 
\begin{displaymath}
	\lim_{\lambda\rightarrow 0}\sin_{\lambda}^{(y)}(t)=\sum_{k=0}^{\infty}(-1)^{k}y^{2k+1}\frac{t^{2k+1}}{(2k+1)!}=\sin(yt).
\end{displaymath}
From \eqref{18}, we note that 
\begin{align}\label{31}
	\sum_{n=0}^{\infty}\mathcal{E}_{n,\lambda}(x+iy)\frac{t^{n}}{n!}&=\frac{2}{e_{\lambda}(t)+1}e^{x}_{\lambda}(t)\cdot e^{iy}_{\lambda}(t)\\
	&=\sum_{l=0}^{\infty}\mathcal{E}_{l,\lambda}(x)\frac{t^{l}}{l!}\sum_{j=0}^{\infty}(iy)_{j,\lambda}\frac{t^{j}}{j!}\nonumber\\
	&=\sum_{n=0}^{\infty}\bigg(\sum_{l=0}^{n}\binom{n}{l}(iy)_{n-l,\lambda}\mathcal{E}_{l,\lambda}(x)\bigg)\frac{t^{n}}{n!}.\nonumber
\end{align}
On the other hand 
\begin{align}\label{32}
	\frac{2}{e_{\lambda}(t)+1}e_{\lambda}^{x+iy}(t)&=\sum_{n=0}^{\infty}\mathcal{E}_{l,\lambda}\frac{t^{l}}{l!}\sum_{j=0}^{\infty}(x+iy)_{j,\lambda}\frac{t^{j}}{j!}\\
	&=\sum_{n=0}^{\infty}\bigg(\sum_{l=0}^{n}\binom{n}{l}(x+iy)_{n-l,\lambda}\mathcal{E}_{l,\lambda}\bigg)\frac{t^{n}}{n!}\nonumber. 
\end{align}
Therefore, by \eqref{31} and \eqref{32}, we obtain the following theorem, 

\begin{theorem}
For $n\ge 0$, we have
\begin{align*}
\mathcal{E}_{n,\lambda}(x+iy)&=\sum_{l=0}^{n}\binom{n}{l}(iy)_{n-l,\lambda}\mathcal{E}_{l,\lambda}(x)\\
&=\sum_{l=0}^{n}\binom{n}{l}(x+iy)_{n-l,\lambda}\mathcal{E}_{l,\lambda}. 
\end{align*}
Also, we have
\begin{align*}
	\mathcal{E}_{n,\lambda}(x-iy)&=\sum_{l=0}^{n}\binom{n}{l}(-1)^{n-l}\langle iy\rangle_{n-l,\lambda}\mathcal{E}_{l,\lambda}(x)\\
	&=\sum_{l=0}^{n}\binom{n}{l}(-1)^{n-l}\langle iy-x\rangle_{n-l,\lambda}\mathcal{E}_{l,\lambda},
\end{align*} 
where $\langle x\rangle_{0,\lambda}=1$, $\langle x\rangle_{n,\lambda}=x(x+\lambda)\cdots(x+\lambda(n-1)),\ (n\ge 1)$. 
\end{theorem}

\noindent By \eqref{29}, we get 
\begin{equation}\label{33}
	e_{\lambda}^{x}(t)\cos_{\lambda}^{(y)}(t)=\sum_{l=0}^{\infty}(x)_{l,\lambda}\frac{t^{l}}{l!}\sum_{m=0}^{\infty}\sum_{k=0}^{[\frac{m}{2}]}\lambda^{m-2k}(-1)^{k}y^{2k}S^{(1)}(m,2k)\frac{t^{m}}{m!}
\end{equation}
\begin{displaymath}
	=\sum_{n=0}^{\infty}\bigg(\sum_{m=0}^{n}\sum_{k=0}^{[\frac{m}{2}]}\binom{n}{m}\lambda^{m-2k}(-1)^{k}y^{2k}S^{(1)}(m,2k)(x)_{n-m,\lambda}\bigg)\frac{t^{n}}{n!}, 
\end{displaymath}
and 
\begin{equation}\label{34}
	e_{\lambda}^{x}(t)\sin_{\lambda}^{(y)}(t)=\sum_{l=0}^{\infty}(x)_{\lambda,l}\frac{t^{l}}{l!}\sum_{m=1}^{\infty}\sum_{k=0}^{[\frac{m-1}{2}]}\lambda^{m-2k-1}(-1)^{k}y^{2k+1}S^{(1)}(m,2k+1)\frac{t^{m}}{m!}
\end{equation}
\begin{displaymath}
	=\sum_{n=1}^{\infty}\bigg(\sum_{m=1}^{n}\sum_{k=0}^{[\frac{m-1}{2}]} \binom{n}{m}\lambda^{m-2k-1}(-1)^{k}y^{2k+1}S^{(1)}(m,2k+1)(x)_{n-m,\lambda}\bigg)\frac{t^{n}}{n!}.
\end{displaymath}
Now, we define the degenerate cosine-polynomials and degenerate sine-polynomials respectively by
\begin{equation}\label{35}
	e_{\lambda}^{x}(t)\cos_{\lambda}^{(y)}(t)=\sum_{k=0}^{\infty}C_{k,\lambda}(x,y)\frac{t^{k}}{k!},
\end{equation}
and 
\begin{equation}\label{36}
	e_{\lambda}^{x}(t)\sin_{\lambda}^{(y)}(t)=\sum_{k=0}^{\infty}S_{k,\lambda}(x,y)\frac{t^{k}}{k!}.
\end{equation}
Note that 
\begin{displaymath}
	\lim_{\lambda\rightarrow 0}C_{k,\lambda}(x,y)=C_{k}(x,y),\quad\lim_{\lambda\rightarrow 0}S_{k,\lambda}(x,y)=S_{k}(x,y),
\end{displaymath}
where $C_{k}(x,y)$ and $S_{k}(x,y)$ are the cosine-polynomials and sine-polynomials of Masijed-Jamei, Beyki and Koepf. \par
~~\\
Therefore, by \eqref{33}-\eqref{36}, we obtain the following theorem. 
\begin{theorem}
	For $n\ge 0$, we have 
	\begin{align*}
		C_{n,\lambda}(x,y)&=\sum_{m=0}^{n}\sum_{k=0}^{[\frac{m}{2}]}\binom{n}{m}\lambda^{m-2k}(-1)^{k}y^{2k}S^{(1)}(m,2k)(x)_{n-m,\lambda}\\
		&= \sum_{k=0}^{[\frac{n}{2}]}\sum_{m=2k}^{n}\binom{n}{m} \lambda^{m-2k}(-1)^{k}y^{2k}S^{(1)}(m,2k)(x)_{n-m,\lambda}.
	\end{align*}
	Also, for $n\in\mathbb{N}$, we have
	\begin{align*}
		S_{n,\lambda}(x,y)&=\sum_{m=1}^{n}\sum_{k=0}^{[\frac{m-1}{2}]}\binom{n}{m}\lambda^{m-2k-1}(-1)^{k}y^{2k+1}S^{(1)}(m,2k+1)(x)_{n-m,\lambda}\\
		&=\sum_{k=0}^{[\frac{n-1}{2}]}\sum_{m=2k+1}^{n} \binom{n}{m}\lambda^{m-2k-1}(-1)^{k}y^{2k+1}S^{(1)}(m,2k+1)(x)_{n-m,\lambda}.
	\end{align*}
and $S_{0,\lambda}(x,y)=0$. 
\end{theorem}
\noindent From \eqref{24}, we note that 
\begin{align}\label{37}
\sum_{n=0}^{\infty}\mathcal{E}_{n,\lambda}^{(c)}(x,y)\frac{t^{n}}{n!}&=\frac{2}{e_{\lambda}(t)+1}e_{\lambda}^{x}(t)\cos_{\lambda}^{(y)}(t)\\
&=\sum_{m=0}^{\infty}\mathcal{E}_{m,\lambda}\frac{t^{m}}{m!}\sum_{l=0}^{\infty}C_{l,\lambda}(x,y)\frac{t^{l}}{l!}\nonumber\\
&=\sum_{n=0}^{\infty}\bigg(\sum_{m=0}^{n}\binom{n}{m}\mathcal{E}_{m,\lambda}C_{n-m,\lambda}(x,y)\bigg)\frac{t^{n}}{n!}\nonumber. 
\end{align}
On the other hand, 
\begin{align}\label{38}
	\frac{2}{e_{\lambda}(t)+1}e_{\lambda}^{x}(t)\cos_{\lambda}^{(y)}(t)&=\sum_{m=0}^{\infty}\mathcal{E}_{m,\lambda}(x)\frac{t^{m}}{m!}\sum_{l=0}^{\infty}\sum_{k=0}^{[\frac{l}{2}]}\lambda^{l-2k}(-1)^{k}y^{2k}S^{(1)}(l,2k)\frac{t^{l}}{l!}\\
 &=\sum_{n=0}^{\infty}\bigg(\sum_{l=0}^{n}\sum_{k=0}^{[\frac{l}{2}]}\binom{n}{l}\lambda^{l-2k}(-1)^{k}y^{2k}S^{(1)}(l,2k)\mathcal{E}_{n-l,\lambda}(x)\bigg)\frac{t^{n}}{n!}\nonumber\\
	 &=\sum_{n=0}^{\infty}\bigg(\sum_{k=0}^{[\frac{n}{2}]}\sum_{l=2k}^{n}\binom{n}{l}\lambda^{l-2k}(-1)^{k}y^{2k}S^{(1)}(l,2k)\mathcal{E}_{n-l,\lambda}(x)\bigg)\frac{t^{n}}{n!}.\nonumber
\end{align}
By \eqref{30}, we get 
\begin{align}\label{39}
\frac{2}{e_{\lambda}(t)+1}e_{\lambda}^{x}(t)\sin_{\lambda}^{(y)}(t)&=\sum_{m=0}^{\infty}\mathcal{E}_{m,\lambda}(x)\frac{t^{m}}{m!}\sum_{l=1}^{n}\sum_{k=0}^{[\frac{l-1}{2}]}(-1)^{k}\lambda^{l-2k-1}y^{2k+1}S^{(1)}(l,2k+1)\frac{t^{l}}{l!}\\
&=\sum_{n=1}^{\infty}\bigg(\sum_{l=1}^{n}\sum_{k=0}^{[\frac{l-1}{2}]}\binom{n}{l}\lambda^{l-2k-1}(-1)^{k}y^{2k+1}S^{(1)}(l,2k+1)\mathcal{E}_{n-l,\lambda}(x)\bigg)\frac{t^{n}}{n!}\nonumber\\
&=\sum_{n=1}^{\infty}\bigg(\sum_{k=0}^{[\frac{n-1}{2}]} \sum_{l=2k+1}^{n}\binom{n}{l}\lambda^{l-2k-1}(-1)^{k}y^{2k+1}S^{(1)}(l,2k+1)\mathcal{E}_{n-l,\lambda}(x)\bigg)\frac{t^{n}}{n!}.\nonumber
\end{align}
Therefore, by \eqref{24}, \eqref{25}, and \eqref{37}-\eqref{39}, we obtain the following theorem. 
\begin{theorem}
	For $n\ge 0$, we have 
	\begin{align*}
		\mathcal{E}_{n,\lambda}^{(c)}(x,y)&=\sum_{k=0}^{n}\binom{n}{k}\mathcal{E}_{k,\lambda}C_{n-k,\lambda}(x,y)\\
		&=\sum_{k=0}^{[\frac{n}{2}]}\sum_{l=2k}^{n}\binom{n}{l}\lambda^{l-2k}(-1)^{k}y^{2k}S^{(1)}(l,2k)\mathcal{E}_{n-l,\lambda}(x).
\end{align*}
Also, for $n\in\mathbb{N}$, we obtain 
\begin{align*}
	\mathcal{E}_{n,\lambda}^{(s)}(x,y)&=\sum_{k=0}^{n}\binom{n}{k}\mathcal{E}_{k,\lambda}S_{n-k,\lambda}(x,y)\\
	&=\sum_{k=0}^{[\frac{n-1}{2}]}\sum_{l=2k+1}^{n}\binom{n}{l}\lambda^{l-2k-1}(-1)^{k}y^{2k+1}S^{(1)}(l,2k+1)\mathcal{E}_{n-l,\lambda}(x). 
\end{align*}
\end{theorem}
By \eqref{24}, we get 
\begin{align}\label{40}	
2e_{\lambda}^{x}(t)\cos_{\lambda}^{(y)}(t)&=\sum_{l=0}^{\infty}\mathcal{E}_{l,\lambda}^{(c)}(x,y)\frac{t^{l}}{l!}(e_{\lambda}(t)+1)\\
&=\sum_{l=0}^{\infty}\mathcal{E}_{l,\lambda}^{(c)}(x,y)\frac{t^{l}}{l!}\sum_{m=0}^{\infty}(1)_{m,\lambda}\frac{t^{m}}{m!}+\sum_{n=0}^{\infty}\mathcal{E}_{n,\lambda}^{(c)}(x,y)\frac{t^{n}}{n!}\nonumber\\
&=\sum_{n=0}^{\infty}\bigg(\sum_{l=0}^{n}\binom{n}{l}(1)_{n-l,\lambda}\mathcal{E}_{l,\lambda}^{(c)}(x,y)+\mathcal{E}_{n,\lambda}^{(c)}(x,y)\bigg)\frac{t^{n}}{n!}.\nonumber
\end{align}
Therefore by comparing the coefficients on both sides of \eqref{35} and \eqref{40}, we obtain the following theorem.
\begin{theorem}
For $n\ge 0$, we have 
\begin{displaymath}
C_{n,\lambda}(x,y)=\frac{1}{2}\bigg(\sum_{l=0}^{n}\binom{n}{l}(1)_{n-l,\lambda}\mathcal{E}_{l,\lambda}^{(c)}(x,y)+\mathcal{E}_{n,\lambda}^{(c)}(x,y)\bigg),
	\end{displaymath}
and
\begin{displaymath}
S_{n,\lambda}(x,y)=\frac{1}{2}\bigg(\sum_{l=0}^{n}\binom{n}{l}(1)_{n-l,\lambda}\mathcal{E}_{l,\lambda}^{(s)}(x,y)+\mathcal{E}_{n,\lambda}^{(s)}(x,y)\bigg).
\end{displaymath}
\end{theorem}
\noindent From \eqref{24}, we have 
\begin{align}\label{41}
\sum_{n=0}^{\infty}\mathcal{E}_{n,\lambda}^{(c)}(x+r,y)\frac{t^{n}}{n!}&=\frac{2}{e_{\lambda}(t)+1}e_{\lambda}^{x+r}(t)\cos_{\lambda}^{(y)}(t)\\
&=\frac{2}{e_{\lambda}(t)+1}e_{\lambda}^{x}(t)\cos_{\lambda}^{(y)}(t) e_{\lambda}^{r}(t)\nonumber\\
&=\sum_{l=0}^{\infty}\mathcal{E}_{l,\lambda}^{(c)}(x,y)\frac{t^{l}}{l!}\sum_{m=0}^{\infty}(r)_{m,\lambda}\frac{t^{m}}{m!}\nonumber\\
&=\sum_{n=0}^{\infty}\bigg(\sum_{l=0}^{n}\binom{n}{l}\mathcal{E}_{l,\lambda}^{(c)}(x,y)(r)_{n-l,\lambda}\bigg)\frac{t^{n}}{n!}.\nonumber
\end{align}
Therefore, by comparing the coefficients on both sides of \eqref{41}, we obtain the following proposition. 
\begin{proposition}
	For $n\ge 0$, we have 
	\begin{displaymath}
		\mathcal{E}_{n,\lambda}^{(c)}(x+r,y)=\sum_{l=0}^{n}\binom{n}{l}\mathcal{E}_{l,\lambda}^{(c)}(x,y)(r)_{n-l,\lambda}, 
	\end{displaymath}
	and 
	\begin{displaymath}
		\mathcal{E}_{n,\lambda}^{(s)}(x+r,y)=\sum_{l=0}^{n}\binom{n}{l}\mathcal{E}_{l,\lambda}^{(s)}(x,y)(r)_{n-l,\lambda}, 
	\end{displaymath}
where $r$ is a fixed real (or complex) number. 
\end{proposition}
\noindent Now, we consider the reflection symmetric identities for the degenerate cosine-Euler polynomials. \\
~~~\
By \eqref{24}, we get 
\begin{align}\label{42}	
\sum_{n=0}^{\infty}\mathcal{E}_{n,\lambda}^{(c)}(1-x,y)\frac{t^{n}}{n!}&=\frac{2}{e_{\lambda}(t)+1}e_{\lambda}^{1-x}(t)\cos_{\lambda}^{(y)}(t)\\
&=\frac{2}{1+e_{\lambda}^{-1}(t)}e_{\lambda}^{-x}(t)\cos_{\lambda}^{(y)}(t)\nonumber\\
&=\frac{2}{e_{-\lambda}(-t)+1}e_{-\lambda}^{x}(-t)\cos_{-\lambda}^{(y)}(-t)\nonumber\\
&=\sum_{n=0}^{\infty}\mathcal{E}_{n,-\lambda}^{(c)}(x,y)\frac{(-1)^{n}t^{n}}{n!}, \nonumber
\end{align}
and 
\begin{align}\label{43}	
\sum_{n=0}^{\infty}\mathcal{E}_{n,\lambda}^{(s)}(1-x,y)\frac{t^{n}}{n!}&=\frac{2}{e_{\lambda}(t)+1}e_{\lambda}^{1-x}(t)\sin_{\lambda}^{(y)}(t)\\
&=\frac{2}{1+e_{\lambda}^{-1}(t)}e_{\lambda}^{-x}(t)\sin_{\lambda}^{(y)}(t)\nonumber\\
&=\frac{2}{e_{-\lambda}(-t)+1}e_{-\lambda}^{x}(-t)\sin_{-\lambda}^{(y)}(-t)\nonumber\\
&=-\sum_{n=0}^{\infty}\mathcal{E}_{n,-\lambda}^{(s)}(x,y)\frac{(-1)^{n}t^{n}}{n!}, \nonumber
\end{align}
Therefore, by \eqref{42} and \eqref{43}, we obtain the following theorem 
\begin{theorem}
	For $n\ge 0$, we have 
	\begin{displaymath}
		\mathcal{E}_{n,\lambda}^{(c)}(1-x,y)=(-1)^{n}\mathcal{E}_{n,-\lambda}^{(c)}(x,y), 
	\end{displaymath}
and 
\begin{displaymath}
		\mathcal{E}_{n,\lambda}^{(s)}(1-x,y)=(-1)^{n+1}\mathcal{E}_{n,-\lambda}^{(s)}(x,y),
\end{displaymath}
\end{theorem}

Now, we observe that 
\begin{align}\label{44}	
\sum_{n=0}^{\infty}\mathcal{E}_{n,\lambda}^{(c)}(x,y)\frac{t^{n}}{n!}&=\frac{2}{e_{\lambda}(t)+1}(e_{\lambda}(t)-1+1)^{x}\cos_{\lambda}^{(y)}(t)\\
&=\frac{2}{e_{\lambda}(t)+1}\sum_{l=0}^{\infty}\binom{x}{l}(e_{\lambda}(t)-1)^{l}\cos_{\lambda}^{(y)}(t)\nonumber\\
&=\frac{2}{e_{\lambda}(t)+1}\cos_{\lambda}^{(y)}(t)\sum_{l=0}^{\infty}(x)_{l}\sum_{k=l}^{\infty}S_{\lambda}^{(2)}(k,l)\frac{t^{k}}{k!}\nonumber\\
&=\sum_{j=0}^{\infty}\mathcal{E}_{j,\lambda}^{(c)}(y)\frac{t^{j}}{j!}\sum_{k=0}^{\infty}\bigg(\sum_{l=0}^{k}(x)_{l}S_{\lambda}^{(2)}(k,l)\bigg)\frac{t^{k}}{k!}\nonumber\\
&=\sum_{n=0}^{\infty}\bigg(\sum_{k=0}^{n}\sum_{l=0}^{k}\binom{n}{k}(x)_{l}S_{\lambda}^{(2)}(k,l)\mathcal{E}_{n-k}^{(c)}(y)\bigg)\frac{t^{n}}{n!}.\nonumber
\end{align}
Therefore, by \eqref{44}, we obtain the following theorem. 
\begin{theorem}
	For $n\ge 0$, we have 
	\begin{displaymath}
		\mathcal{E}_{n,\lambda}^{(c)}(x,y)=\sum_{k=0}^{n}\sum_{l=0}^{k}\binom{n}{l}(x)_{l}S_{\lambda}^{(2)}(k,l)\mathcal{E}_{n-k,\lambda}^{(c)}(y). 
	\end{displaymath}
Also, for $n\in\mathbb{N}$, we have 
\begin{displaymath}
	\mathcal{E}_{n,\lambda}^{(s)}(x,y)=\sum_{k=0}^{n}\sum_{l=0}^{k}\binom{n}{k}(x)_{l}S_{\lambda}^{(2)}(k,l)\mathcal{E}_{n-k,\lambda}^{(s)}(y). 
\end{displaymath}
\end{theorem}

\section{Degenerate Bernoulli polynomials of complex variable} 

In this section, we will consider the degenerate Bernoulli polynomials of complex variable and, by treating the real and imaginary parts separately, introduce the degenerate cosine-Bernoulli polynomials and degenerate sine-Bernoulli polynomials.

\noindent From \eqref{4}, we have 
\begin{equation}\label{45}
	\frac{t}{e_{\lambda}(t)-1}e_{\lambda}^{x+iy}(t)=\sum_{n=0}^{\infty}\beta_{n,\lambda}(x+iy)\frac{t^{n}}{n!},
\end{equation}
and 
\begin{equation}\label{46}
	\frac{t}{e_{\lambda}(t)-1}e_{\lambda}^{x-iy}=\sum_{n=0}^{\infty}\beta_{n,\lambda}(x-iy)\frac{t^{n}}{n!}.
\end{equation}
Thus, by \eqref{45} and \eqref{46}, we get 
\begin{equation}\label{47}
	\sum_{n=0}^{\infty}\big(\beta_{n,\lambda}(x+iy)+\beta_{n,\lambda}(x-iy)\big)\frac{t^{n}}{n!}=2\frac{t}{e_{\lambda}(t)-1}e_{\lambda}^{x}(t)\cos_{\lambda}^{(y)}(t),
\end{equation}
and 
\begin{equation}\label{48}
	\sum_{n=0}^{\infty}\big(\beta_{n,\lambda}(x+iy)-\beta_{n,\lambda}(x-iy)\big)\frac{t^{n}}{n!}=2i\frac{t}{e_{\lambda}(t)-1}e_{\lambda}^{x}(t)\sin_{\lambda}^{(y)}(t).
\end{equation}
In view of \eqref{24} and \eqref{25}, we define the degenerate cosine-Bernoulli polynomials and degenerate sine-Bernoulli polynomials respectively by
\begin{equation}\label{49}
	\frac{t}{e_{\lambda}(t)-1}e^{x}_{\lambda}(t)\cos_{\lambda}^{(y)}(t)=\sum_{n=0}^{\infty}\beta_{n,\lambda}^{(c)}(x,y)\frac{t^{n}}{n!},
\end{equation}
and
\begin{equation}\label{50}
	\frac{t}{e_{\lambda}(t)-1}e^{x}_{\lambda}(t)\sin_{\lambda}^{(y)}(t)=\sum_{n=0}^{\infty}\beta_{n,\lambda}^{(s)}(x,y)\frac{t^{n}}{n!}.
\end{equation}
Note that $\beta_{0,\lambda}^{(s)}(x,y)=0$.\\ 
From \eqref{47}-\eqref{50}, we have 
\begin{equation}\label{51}
	\beta_{n,\lambda}^{(c)}(x,y)=\frac{\beta_{n,\lambda}(x+iy)+\beta_{n,\lambda}(x-iy)}{2},
\end{equation}
and 
\begin{equation}\label{52}
	\beta_{n,\lambda}^{(s)}(x,y)=\frac{\beta_{n,\lambda}(x+iy)-\beta_{n,\lambda}(x-iy)}{2i},\quad (n\ge 0).
\end{equation}
Note that 
\begin{displaymath}
	\lim_{\lambda\rightarrow 0}\beta_{n,\lambda}^{(c)}(x,y)=B_{n}^{(c)}(x,y),\quad\lim_{\lambda\rightarrow 0}\beta_{n,\lambda}^{(s)}(x,y)=B_{n}^{(s)}(x,y),
\end{displaymath}
where $B_{n}^{(c)}(x,y),\ B_{n}^{(s)}(x,y)$ are cosine-Bernoulli polynomials, and sine-Bernoulli polynomials (see [12,16]). ~~\\
~~\\
By \eqref{49}, we get 
\begin{align}\label{53}	
\sum_{n=0}^{\infty}\beta_{n,\lambda}^{(c)}(x,y)\frac{t^{n}}{n!}&=\frac{t}{e_{\lambda}(t)-1}e_{\lambda}^{x}(t)\cos_{\lambda}^{(y)}(t)\\
&=\sum_{l=0}^{\infty}\beta_{l,\lambda}\frac{t^{l}}{l!}\sum_{m=0}^{\infty}C_{m,\lambda}(x,y)\frac{t^{m}}{m!}\nonumber\\
&=\sum_{n=0}^{\infty}\bigg(\sum_{l=0}^{n}\binom{n}{l}\beta_{l,\lambda}C_{n-l,\lambda}(x,y)\bigg)\frac{t^{n}}{n!}.\nonumber
\end{align}
On the other hand, 
\begin{align}\label{54}	
\frac{t}{e_{\lambda}(t)-1}&e_{\lambda}^{x}(t)\cos_{\lambda}^{(y)}(t)\\
&=\sum_{m=0}^{\infty}\beta_{m,\lambda}(x)\frac{t^{m}}{m!}\sum_{l=0}^{n}\sum_{k=0}^{[\frac{l}{2}]}\lambda^{l-2k}(-1)^{k}y^{2k}S^{(1)}(l,2k)\frac{t^{l}}{l!}\nonumber\\
&=\sum_{n=0}^{\infty}\bigg(\sum_{l=0}^{n}\sum_{k=0}^{[\frac{l}{2}]}\binom{n}{l}\lambda^{l-2k}(-1)^{k}y^{2k}S^{(1)}(l,2k)\beta_{n-l,\lambda}(x)\bigg)\frac{t^{n}}{n!}\nonumber \\
&=\sum_{n=0}^{\infty}\bigg(\sum_{k=0}^{[\frac{n}{2}]}\sum_{l=2k}^{n} \binom{n}{l}\lambda^{l-2k}(-1)^{k}y^{2k}S^{(1)}(l,2k)\beta_{n-l,\lambda}(x)\bigg)\frac{t^{n}}{n!}.\nonumber
\end{align}
Therefore, by \eqref{53} and \eqref{54}, we obtain the following theorem.
\begin{theorem}
	For $n\ge 0$, we have 
	\begin{align*}
		\beta_{n,\lambda}^{(c)}(x,y)&=\sum_{k=0}^{n}\binom{n}{k}\beta_{k,\lambda}C_{n-k,\lambda}(x,y)\\
		&=\sum_{k=0}^{[\frac{n}{2}]}\sum_{l=2k}^{n}\binom{n}{l}\lambda^{l-2k}(-1)^{k}y^{2k}S^{(1)}(l,2k)\beta_{n-l,\lambda}(x).
	\end{align*}
	Also, for $n\in\mathbb{N}$, we have 
\begin{align*}
		\beta_{n,\lambda}^{(s)}(x,y)&=\sum_{k=0}^{n}\binom{n}{k}\beta_{k,\lambda}S_{n-k,\lambda}(x,y)\\
		&=\sum_{k=0}^{[\frac{n-1}{2}]}\sum_{l=2k+1}^{n}\binom{n}{l}\lambda^{l-2k-1}(-1)^{k}y^{2k+1}S^{(1)}(l,2k+1)\beta_{n-l,\lambda}(x).
	\end{align*}	
and 
\begin{displaymath}
	\beta_{0,\lambda}^{(s)}(x,y)=0.
\end{displaymath}
\end{theorem}
\noindent From \eqref{49}, we have 
\begin{align} \label{55}	
\sum_{n=0}^{\infty}\beta_{n,\lambda}^{(c)}(1-x,y)\frac{t^{n}}{n!}
&=\frac{t}{1-e_{\lambda}^{-1}(t)} e_{\lambda}^{-x}(t)\cos_{\lambda}^{(y)}(t)\\
&=\frac{-t}{e_{-\lambda}(-t)-1}e_{-\lambda}^{x}(-t)\cos_{-\lambda}^{(y)}(-t)\nonumber\\
&= \sum_{n=0}^{\infty}\beta_{n,-\lambda}^{(c)}(x,y)\frac{(-1)^{n}}{n!}t^{n}.\nonumber
\end{align}
Therefore, by \eqref{55}, we obtain the following theorem.
\begin{theorem}
	For $n\ge 0$, we have 
	\begin{displaymath}
		\beta_{n,\lambda}^{(c)}(1-x,y)=(-1)^{n}\beta_{n,-\lambda}^{(c)}(x,y), 
	\end{displaymath}
	and 
	\begin{displaymath}
		\beta_{n,\lambda}^{(s)}(1-x,y)= (-1)^{n+1}\beta_{n,-\lambda}^{(s)}(x,y).
	\end{displaymath}
\end{theorem}
\noindent By \eqref{49}, we easily get 
\begin{align}\label{56}	
\sum_{n=0}^{\infty}\beta_{n,\lambda}^{(c)}(x+r,y)\frac{t^{n}}{n!}&=\frac{t}{e_{\lambda}(t)-1}e_{\lambda}^{x+r}(t)\cos_{\lambda}^{(y)}(t)\\
&=\frac{t}{e_{\lambda}(t)-1}e_{\lambda}^{x}(t)\cos_{\lambda}^{(y)}(t)e_{\lambda}^{r}(t)\nonumber\\
&=\sum_{l=0}^{\infty}\beta_{l,\lambda}^{(c)}(x,y)\frac{t^{l}}{l!}\sum_{m=0}^{\infty}(r)_{m,\lambda}\frac{t^{m}}{m!}\nonumber\\
&=\sum_{n=0}^{\infty}\bigg(\sum_{l=0}^{n}\binom{n}{l}\beta_{l,\lambda}^{(c)}(x,y)(r)_{n-l,\lambda}\bigg)\frac{t^{n}}{n!}\nonumber.	
\end{align}
By comparing the coefficients on both sides of \eqref{56}, we get 
\begin{equation}\label{57}
	\beta_{n\lambda}^{(c)}(x+r,y)=\sum_{l=0}^{n}\binom{n}{l}\beta_{l,\lambda}^{(c)}(x,y)(r)_{n-l,\lambda},
\end{equation}
and 
\begin{equation}\label{58}
	\beta_{n,\lambda}^{(s)}(x+r,y)=\sum_{l=0}^{n}\binom{n}{l}\beta_{l,\lambda}^{(s)}(x,y)(r)_{n-l,\lambda},
\end{equation}
where $r$ is a fixed real (or complex) number. \\
~~\\
From \eqref{49}, we note that 
\begin{align}\label{59}	
te_{\lambda}^{x}(t)\cos_{\lambda}^{(y)}(t)&=\sum_{l=0}^{\infty}\beta_{l,\lambda}^{(c)}(x,y)\frac{t^{l}}{l!}(e_{\lambda}(t)-1)\\
&= \sum_{l=0}^{\infty}\beta_{l,\lambda}^{(c)}(x,y)\frac{t^{l}}{l!}\sum_{m=0}^{\infty}(1)_{m,\lambda}\frac{t^{m}}{m!}-\sum_{n=0}^{\infty}\beta_{n,\lambda}^{(c)}(x,y)\frac{t^{n}}{n!}\nonumber \\
&=\sum_{n=0}^{\infty}\bigg(\sum_{l=0}^{n}\binom{n}{l}\beta_{l,\lambda}^{(c)}(x,y)(1)_{n-l,\lambda}-\beta_{n,\lambda}^{(c)}(x,y)\bigg)\frac{t^{n}}{n!}\nonumber \\
&=\sum_{n=1}^{\infty}\bigg(\beta_{n,\lambda}^{(c)}(x+1,y)-\beta_{n,\lambda}^{(c)}(x,y)\bigg)\frac{t^{n}}{n!}\nonumber\\
&=\sum_{n=0}^{\infty}\bigg(\frac{\beta_{n+1,\lambda}^{(c)}(x+1,y)-\beta_{n+1,\lambda}^{(c)}(x,y)}{n+1}\bigg)\frac{t^{n+1}}{n!}.\nonumber
\end{align}
By \eqref{59}, we get 
\begin{equation}\label{60}
\sum_{n=0}^{\infty}\bigg(\frac{\beta_{n+1,\lambda}^{(c)}(x+1,y)-\beta_{n+1,\lambda}^{(c)}(x,y)}{n+1}\bigg)\frac{t^{n}}{n!}=e_{\lambda}^{x}(t)\cos_{\lambda}^{(y)}(t)=\sum_{n=0}^{\infty}C_{n,\lambda}(x,y)\frac{t^{n}}{n!}.
\end{equation}
Therefore, by comparing the coefficients on both sides of \eqref{60}, we obtain the following theorem.
\begin{theorem}
	For $n\ge 0$, we have 
	\begin{displaymath}
		C_{n,\lambda}(x,y)=\frac{1}{n+1}\big\{\beta_{n+1,\lambda}^{(c)}(x+1,y)-\beta_{n+1,\lambda}^{(c)}(x,y)\big\},
	\end{displaymath}
	and 
	\begin{displaymath}
		S_{n,\lambda}(x,y)=\frac{1}{n+1}\big\{\beta_{n+1,\lambda}^{(s)}(x+1,y)-\beta_{n+1,\lambda}^{(s)}(x,y)\big\}.
	\end{displaymath}
\end{theorem}
\begin{corollary}
For $n\ge 1$, we have 
	\begin{displaymath}
		C_{n,\lambda}(x,y)=\frac{1}{n+1}\sum_{l=0}^{n}\binom{n+1}{l}\beta_{l,\lambda}^{(c)}(x,y)(1)_{n+1-l,\lambda}, 
	\end{displaymath}
	and 
	\begin{displaymath}
		S_{n,\lambda}(x,y)= \frac{1}{n+1}\sum_{l=0}^{n}\binom{n+1}{l}\beta_{l,\lambda}^{(s)}(x,y)(1)_{n+1-l,\lambda}.
	\end{displaymath}
\end{corollary}
\noindent When $x=0$, let $\beta_{n,\lambda}^{(c)}(0,y)=\beta_{n,\lambda}^{(c)}(y),\ \beta_{n,\lambda}^{(s)}(0,y)=\beta_{n,\lambda}^{(s)}(y)$, $\mathcal{E}_{n,\lambda}^{(c)}(0,y)=\mathcal{E}_{n,\lambda}^{(c)}(y)$, and\\
$\mathcal{E}_{n,\lambda}^{(s)}(0,y)=\mathcal{E}_{n,\lambda}^{(s)}(y)$. \\
~~~\\
For $n\ge 0$, we have 
\begin{equation}\label{61}
	\beta_{n,\lambda}^{(c)}(y)=\sum_{k=0}^{[\frac{n}{2}]}\sum_{l=2k}^{n}\binom{n}{l}\lambda^{l-2k}(-1)^{k}y^{2k}S^{(1)}(l,2k)\beta_{n-l,\lambda}. 
\end{equation}
Also, for $n\in\mathbb{N}$, we get 
\begin{equation}\label{62}
	\beta_{n,\lambda}^{(s)}(y)=\sum_{k=0}^{[\frac{n-1}{2}]}\sum_{l=2k+1}^{n}\binom{n}{l}\lambda^{l-2k-1}(-1)^{k}y^{2k+1}S^{(1)}(l,2k+1)\beta_{n-l,\lambda}.
\end{equation}

\noindent By \eqref{49}, we get 
\begin{align}\label{63}
\sum_{n=0}^{\infty}\beta_{n,\lambda}^{(c)}(x,y)\frac{t^{n}}{n!}&=\frac{t}{e_{\lambda}(t)-1}\cos_{\lambda}^{(y)}(t)\big(e_{\lambda}(t)-1+1\big)^{x}\\
&=\sum_{m=0}^{\infty}\beta_{m,\lambda}^{(c)}(y)\frac{t^{m}}{m!}\sum_{l=0}^{\infty}(x)_{l}\sum_{k=l}^{\infty}S_{\lambda}^{(2)}(k,l)\frac{t^{k}}{k!}\nonumber\\
&=\sum_{m=0}^{\infty}\beta_{m,\lambda}^{(c)}(y)\frac{t^{m}}{m!}\sum_{k=0}^{\infty}\sum_{l=0}^{k}(x)_{l}S_{\lambda}^{(2)}(k,l)\frac{t^{k}}{k!}\nonumber\\
&=\sum_{n=0}^{\infty}\bigg(\sum_{k=0}^{n}\sum_{l=0}^{k}\binom{n}{k}(x)_{l}S_{\lambda}^{(2)}(k,l)\beta_{n-k,\lambda}^{(c)}(y)\bigg)\frac{t^{n}}{n!}.\nonumber
\end{align}
Comparing the coefficients on both sides of \eqref{63}, we have 
\begin{displaymath}
	\beta_{n,\lambda}^{(c)}(x,y)=\sum_{k=0}^{n}\sum_{l=0}^{k}\binom{n}{k}(x)_{l}S_{\lambda}^{(2)}(k,l)\beta_{n-k,\lambda}^{(c)}(y).
\end{displaymath}
Also, for $n\in\mathbb{N}$, we get 
\begin{displaymath}
	\beta_{n,\lambda}^{(s)}(x,y)=\sum_{k=0}^{n}\sum_{l=0}^{k}\binom{n}{k}(x)_{l}S_{\lambda}^{(2)}(k,l)\beta_{n-k,\lambda}^{(s)},
\end{displaymath}
and 
\begin{displaymath}
	\beta_{0,\lambda}^{(s)}(x,y)=0. 
\end{displaymath}

\section{Conclusions}
In [15], the authors introduced the so called the new type Euler polynomials by means of generating functions (see \eqref{9}, \eqref{10}) and deduced several properties and identities for these polynomials. Hac\`ene Belbachir, the reviewer of the paper [15], asked the following question in Mathematical Reviews (MR3808565) of the American Mathematical Society:
Is it possible to obtain their results by considering the classical Euler polynomials of complex variable $z$, and treating the real part and the imaginary part separately?\\
\indent Our result gives an affirmative answer to the question (see \eqref{16}). In this paper, we considered the degenerate Euler and Bernoulli polynomials of complex variable and, by treating the real and imaginary parts separately, were able to introduce degenerate cosine-Euler polynomials, degenerate sine-Euler polynomials, degenerate cosine-Bernoulli polynomials, and degenerate sine-Bernoulli polynomials. They are degenerate versions of the new type Euler polynomials studied by  Masjed-Jamei, Beyki and Koepf [15] and of the 'new type Bernoulli polynomials.'\\ 
\indent In Section 2, the degenerate cosine-polynomials and degenerate sine-polynomials were introduced and their explicit expressions were derived. The degenerate cosine-Euler polynomials and degenerate sine-Euler polynomials were expressed in terms of degenerate cosine-polynomials and degenerate sine-polynomials and vice versa. Further, some reflection identities were found for the degenerate cosine-Euler polynomials and degenerate sine-Euler polynomials. In Section 3, the degenerate cosine-Bernoulli polynomials and degenerate sine-Bernoulli polynomials were introduced. They were expressed in terms of degenerate cosine-polynomials and degenerate sine-polynomials and vice versa. Reflection symmetries were deduced for the degenerate cosine-Bernoulli polynomials and degenrate sine-Bernoulli polynomials.
Further, some expressions involving the degenerate Stirling numbers of the second kind were derived  for them.\\
\indent It was Carlitz [1,2] who initiated the study of degenerate versions of some special polynomials, namely the degenerate Bernoulli and Euler polynomials. Studying degenerate versions of some special polynomials and numbers have turned out to be very fruitful and promising (see [3,5-11,13-14,19] and references therein). In fact, this idea of considering degenerate versions of some special polynomials are not limited just to polynomials but can be extended even to transcendental functions like gamma functions [8].

\end{document}